\documentclass[11pt,smallextended]{article}
\usepackage[latin1]{inputenc}
\usepackage[english]{babel}
\usepackage{geometry}
\usepackage{amsmath}
\usepackage{amssymb}
\usepackage{oldgerm}
\usepackage{amscd}
\usepackage{amsthm}
\usepackage{rotating}
\usepackage[all]{xy}
\usepackage{hyperref}

\title{Integer-valued polynomials over matrices and divided differences}

\date{\today}

\author{Giulio Peruginelli\footnote{Institut f\"ur Analysis und Comput. Number Theory, Technische Univ., Steyrergasse 30, A-8010 Graz, Austria. E-mail: peruginelli@math.tugraz.at.}}

\theoremstyle{definition}\newtheorem{Def}{Definition}[section]
\newtheorem{Ex}{Example}[section]
\newtheorem{Rem}{Remark}[section]

\theoremstyle{plain}\newtheorem{Prop}{Proposition}[section]
\newtheorem{Theorem}{Theorem}[section]
\newtheorem{Lemma}{Lemma}[section]
\newtheorem{Cor}{Corollary}[section]
\newcommand{\Proof}{\noindent{\bfseries Proof : }} 
\newcommand{\Q}{\mathbb{Q}}
\newcommand{\N}{\mathbb{N}}
\newcommand{\Z}{\mathbb{Z}}
\newcommand{\IZ}{\textnormal{Int($\Z$)}}

\begin{document}
\maketitle
\begin{abstract}
\noindent 
Let $D$ be an integrally closed domain with quotient field $K$ and $n$ a positive integer. We give a characterization of the polynomials in $K[X]$ which are integer-valued over the set of matrices $M_n(D)$ in terms of their divided differences. A necessary and sufficient condition on $f\in K[X]$ to be integer-valued over $M_n(D)$ is that, for each $k$ less than $n$, the $k$-th divided difference of $f$ is integral-valued on every subset of the roots of any monic polynomial over $D$ of degree $n$. If in addition the intersection of the maximal ideals of finite index is $(0)$ then it is sufficient to check the above conditions on subsets of the roots of monic irreducible polynomials of degree $n$, that is, conjugate integral elements of degree $n$ over $D$.\footnote{The only difference from this version and the official one (\href{http://dx.doi.org/10.1007/s00605-013-0519-9}{http://dx.doi.org/10.1007/s00605-013-0519-9}) is that the assumption "$D$ has zero Jacobson radical" in (\ref{Mnirr}), (\ref{IntMnIntMnpzerJAc}) and Thm \ref{MTINTMNGLOBALE} has to be replaced with "the intersection of the maximal ideals of finite index is $(0)$". }
\end{abstract}
\vskip0.3cm


\vskip0.5cm
\section{Introduction}

Let $D$ be a domain with quotient field $K$. We denote by $M_n(D)$ the $D$-algebra of $n\times n$ matrices  with entries in $D$. We consider the ring of polynomials in $K[X]$ which are integer-valued over $M_n(D)$ (see \cite{Frisch1}, for a general reference), namely:

$$\textnormal{Int}(M_n(D))\doteqdot\{f\in K[X]\,|\,f(M)\in M_n(D),\,\forall M\in M_n(D)\}.$$

In \cite{EvFaJoh} the authors consider the overring of integer-valued polynomials which are integer-valued over the subset $T_n(D)\subset M_n(D)$ of triangular matrices over $D$, namely the ring $\textnormal{Int}(T_n(D))\doteqdot\{f\in K[X]\,|\,f(T)\in M_n(D),\,\forall T\in T_n(D)\}$. They characterize a polynomial $f(X)$ of this ring in terms of its divided differences, proving that up to the order $n-1$ the $k$-th divided difference $\Phi^k(f)\in K[X_1,\ldots,X_{k+1}]$ is integer-valued over $D^{k+1}$, that is $\Phi^k(f)(D^{k+1})\subset D$. See section \ref{divdiff} of this paper for the definition of divided differences of a polynomial.

We give here an analogous characterization for a polynomial $f(X)$ in $\textnormal{Int}(M_n(D))$. It turns out that, for every $0\leq k<n$, $\Phi^k(f)$ is integral-valued on every subset of $k+1$ elements of the roots of any monic polynomial $p\in D[X]$ of degree $n$. More precisely, given such a polynomial $p(X)$, let $\alpha_1,\ldots,\alpha_n$ be its roots in a fixed algebraic closure of $K$, counted with multiplicity.  Then $\Phi^k(f)(\alpha_1,\ldots,\alpha_{k+1})$ is integral over $D$, for every possible choice of the subset $\{\alpha_i\}_{i=1}^{k+1}$ of the roots of $p(X)$. This property of $\Phi^k(f)$ holds for every polynomial $p(X)$ chosen as before. Notice that the $\alpha_i$'s are integral over $D$ since $p(X)$ is monic (they are not conjugate over $D$, if $p(X)$ is reducible). These properties about the divided differences of $f(X)$ are sufficient and necessary for $f(X)$ to belong to $\textnormal{Int}(M_n(D))$. This characterization generalizes the previous one about integer-valued polynomials over triangular matrices. In fact, for such a polynomial $f(X)$ and for all the relevant $k$'s, $\Phi^k(f)$ is integer-valued on every subset of $k+1$ roots of every monic polynomial over $D$ of degree $n$ which is totally split over $D$.

Here is an overview of the paper. In the second section we give a characterization of the ring of  polynomials which are integer-valued over the set of matrices with prescribed characteristic polynomial: it turns out that this ring has a quite simple algebraic structure, being equal to a pullback of $K[X]$ (Lemma \ref{crit0} and Remark \ref{pullb}). In the third section, under the assumption that $D$ is integrally closed, we prove that an ideal in $D[X]$ generated by a non-zero element of $D$ and a non-constant monic polynomial is contracted with respect to any polynomial ring extension $D'[X]$, under a suitable assumption. In the fourth section we recall the definition of divided differences and in Theorem \ref{MTINTMNGLOBALE} we use the previous results to characterize integer-valued polynomials over matrices in terms of their divided differences. In the fifth section we give some results about integer-valued polynomials over triangular matrices. In the sixth section we give some results about the so-called polynomial closure of set of matrices having prescribed characteristic polynomial. Finally, in the last section we give some remarks which follows from our main Theorem \ref{MTINTMNGLOBALE}.

\vskip0.5cm

\section{Integer-valued polynomials over matrices with prescribed characteristic polynomial}
\vskip0.4cm

Given a commutative ring $R$ and a matrix $M\in M_n(R)$, we denote by $p_M\in R[X]$ the characteristic polynomial of $M$ and by $N_{R[X]}(M)$ the ideal of polynomials $g(X)$ in $R[X]$ such that $g(M)=0$, also called the null ideal of $M$ over $R$.

As in \cite{EvFaJoh}, given a subset $S$ of $M_n(D)$ we consider the ring of polynomials which are integer-valued over $S$:
$${\rm Int}(S,M_n(D))\doteqdot\{f\in K[X]\,|\,f(M)\in M_n(D),\forall M\in S\}.$$
Notice that we have the containment $\textnormal{Int}(M_n(D))\subseteq{\rm Int}(S,M_n(D))$.

We introduce the following definition:

\begin{Def}\label{Mnp}
For a given positive integer $n$, we denote by $\mathcal P_n$ the set of monic polynomials of degree $n$ in $D[X]$. Given a subset $P\subseteq \mathcal P_n$, we denote by 
$$M_n^P(D)\doteqdot\{M\in M_n(D)\,|\,p_M(X)\in P\}$$
the set of matrices in $M_n(D)$ whose characteristic polynomial $p_M(X)$ is in $P$. To ease the notation, given a polynomial $p\in \mathcal P_n$, we set $M_n^p(D)=M_n^{\{p\}}(D)$.
\end{Def}

\noindent Notice that $M_n^p(D)$ is always non-empty, since the companion matrix $C_p$ of $p(X)$ belongs to $M_n^p(D)$. Given $p\in \mathcal P_n$, we begin to study the corresponding ring of integer-valued polynomials $\textnormal{Int}(M_n^p(D),M_n(D))$.

The following Lemma is easy to prove.

\begin{Lemma}\label{pNM}
Let $R$ be a commutative ring and $p\in R[X]$ a monic polynomial of degree $n$. Then
  $$\bigcap_{M\in M_n^p(R)}N_{R[X]}(M)=(p(X)).$$
\end{Lemma}
\Proof By Cayley-Hamilton theorem we have that $p(X)$ annihilates every matrix in $M_n^p(R)$. Conversely, we use the fact that over any commutative ring $R$ the null ideal of the companion matrix of a polynomial $p(X)$ is the principal ideal generated by $p(X)$ (see \cite{Frisch0}). Hence, the above intersection is contained in $N_{R[X]}(C_p)=(p(X))$, so we are done. $\Box$\vskip0.5cm
The next lemma follows immediately from Lemma \ref{pNM}. It gives a characterization of the ring of integer-valued polynomials over the set of matrices with a prescribed characteristic polynomial.
\vskip0.4cm
\begin{Lemma}\label{crit0}
 Let $p\in D[X]$ be a monic polynomial of degree $n$ and $f(X)=g(X)/d\in K[X]$, $g\in D[X]$, $d\in D\setminus\{0\}$. Then the following are equivalent:
\begin{itemize}
 \item[1)] $f\in\textnormal{Int}(M_n^p(D),M_n(D))$.
 \item[2)] $g(X)$ is divisible by $p(X)$ modulo $dD[X]$.
 \item[3)] the remainder of the division of $f(X)$ by $p(X)$ in $K[X]$ is in $D[X]$, that is: $$f(X)=r(X)+p(X)q(X),$$
 for some $q\in K[X]$, $p\in D[X]$ with $\deg(p)<n$.
 \item[4)] $f\in\textnormal{Int}(\{C_p\},M_n(D))$, where $C_p$ is the companion matrix of $p(X)$.
\end{itemize} 
\end{Lemma}
\Proof Given $f(X)=g(X)/d\in K[X]$, it is easy to see that (see also \cite{Frisch0}) we have 
\begin{equation}\label{fMgN}
f(M)\in M_n(D) \textnormal{ if and only if } \overline{g}\in N_{(D/dD)[X]}(\overline{M}),
\end{equation}
where $N_{(D/dD)[X]}(\overline{M})$ is the null ideal over $(D/dD)[X]$ of the matrix $\overline{M}$ obtained by reducing the entries of $M$ modulo the ideal $dD$ and $\overline{g}$ is the polynomial obtained by reducing the coefficients of $g$ modulo $dD$. Since $\textnormal{Int}(M_n^p(D),M_n(D))=\bigcap_{M\in M_n^p(D)}\textnormal{Int}(\{M\},M_n(D))$, we get
\begin{equation}\label{MDp}
\textnormal{Int}(M_n^p(D),M_n(D))=\{f(X)=g(X)/d\in K[X]\;|\;\overline{g}\in\bigcap_{M\in M_n^p(D)}N_{(D/dD)[X]}(\overline{M})\}
\end{equation}
(notice that given $g(X)/d\in K[X]$, the reduction on the right is modulo $d$). By Lemma \ref{pNM} we get that 1)$\Leftrightarrow$ 2)\footnote{This proof can be shortened using directly the above result that over any commutative ring $R$, the null ideal of the companion matrix of a polynomial $p(X)$ is the principal ideal generated by $p(X)$.}. 

Since the polynomial $p(X)$ is monic, by the Euclidean division algorithm condition 2) is equivalent to the condition that $g(X)$ belongs to the ideal of $D[X]$ generated by $p(X)$ and $d$. Clearly, this is equivalent to condition 3).

Finally, by the aforementioned result about null ideals of companion matrices, we have
$$\textnormal{Int}(M_n^p(D),M_n(D))=\textnormal{Int}(\{C_p\},M_n(D)).$$
$\Box$
\vskip0.3cm

\begin{Rem}\label{pullb} We notice that condition 3) says that
$$\textnormal{Int}(M_n^p(D),M_n(D))=D[X]+p(X)\cdot K[X],$$
where the latter ring is a pullback of $K[X]$. In particular, a polynomial $f\in K[X]$ of degree less than $n$ is integer-valued over $M_n^p(D)$ if and only if $f(X)$ is in $D[X]$.
\end{Rem}

It is easy to generalize Lemma \ref{crit0} to integer-valued polynomials over $M_n^P(D)$, for a given subset $P\subseteq\mathcal{P}_n$. Since $M_n^P(D)=\bigcup_{p\in P}M_n^p(D)$, we have 
$$\textnormal{Int}(M_n^P(D),M_n(D))=\bigcap_{p\in P}\textnormal{Int}(M_n^p(D),M_n(D)).$$
The following Proposition follows by Lemma \ref{crit0} and the previous representation of the ring $\textnormal{Int}(M_n^P(D),M_n(D))$ as an intersection of the overrings $\textnormal{Int}(M_n^p(D),M_n(D))$.
 
\begin{Prop}\label{criterio}
Let $P$ be a set of monic polynomials in $D[X]$ of degree $n$.  Let $f(X)=g(X)/d\in K[X]$. Then
$$f\in\textnormal{Int}(M_n^P(D),M_n(D))\Leftrightarrow g(X) \textnormal{ is divisible modulo }dD[X] \textnormal{ by all }p\in P.$$
\end{Prop}
\vskip0.3cm
\begin{Rem} By \cite[Lemma 3.4]{Frisch0} a polynomial $f(X)=g(X)/d\in K[X]$ is integer-valued over $M_n(D)$ if and only if $g(X)$ is divisible modulo $dD[X]$ by all monic polynomials $p\in D[X]$ of degree $n$. Then by Lemma \ref{crit0} we have this representation of the ring $\textnormal{Int}(M_n(D))$:
\begin{equation}\label{IntMnIntMnpnozerJAc}
\textnormal{Int}(M_n(D))=\bigcap_{p\in \mathcal P_n}\textnormal{Int}(M_n^p(D),M_n(D))
\end{equation}
Notice that this is Proposition \ref{criterio} for $P=\mathcal P_n$.

Suppose now that the intersection of the maximal ideals of finite index is $(0)$. By Proposition 6.2 of \cite{Frisch1}\footnote{We notice that in \cite[Proposition 6.2]{Frisch1} the assumption "$D$ has zero Jacobson radical" has to be replaced with "the intersection of the maximal ideals of finite index is $(0)$". Indeed, that theorem is clearly false if for example $K$ is algebraically closed.}  we can just check the previous condition on irreducible polynomials; namely, a polynomial $f(X)=g(X)/d\in K[X]$ is integer-valued over $M_n(D)$ if and only if $g(X)$ is divisible modulo $dD[X]$ by all monic irreducible polynomials $p\in D[X]$ of degree $n$. For a given integer $n$, we denote by $\mathcal P_n^{\textnormal{irr}}$ the set of monic and irreducible polynomials in $D[X]$ of degree $n$. We consider the set of matrices $M_n^{\textnormal{irr}}(D)\subsetneq M_n(D)$  with irreducible characteristic polynomial. Then the last result shows that
\begin{equation}\label{Mnirr}
\textnormal{Int}(M_n(D))=\textnormal{Int}(M_n^{\textnormal{irr}}(D),M_n(D))
\end{equation}
that is, $M_n^{\textnormal{irr}}(D)$ is polynomially dense in $M_n(D)$ (see \cite{EvFaJoh} for the definition, which we recall later in section \ref{polclomatrices}). Hence, in this case, since $M_n^{\textnormal{irr}}(D)=\bigcup_{p\in \mathcal P_n^{\textnormal{irr}}}M_n^p(D)$
we have:
\begin{equation}\label{IntMnIntMnpzerJAc}
\textnormal{Int}(M_n(D))=\bigcap_{p\in \mathcal P_n^{\textnormal{irr}}}\textnormal{Int}(M_n^p(D),M_n(D)).
\end{equation}
\end{Rem}

\vskip0.3cm
\section{A contracted ideal in a polynomial ring extension}
\vskip0.3cm
The next lemma says the following. Let $p(X)$ a monic polynomial over $D$. If a polynomial $g\in D[X]$ is divisible by $p(X)$ modulo an element $d$ of $D$ over a ring extension $D'$ of $D$ (under a suitable hypothesis), then the same modular divisibility condition for $g(X)$ already holds over $D$ itself. Notice that this is exactly condition 2) of Lemma \ref{crit0}.
\vskip0.3cm

\begin{Lemma}\label{gdpK}
Let $D\subset D'$ be an extension of domains with quotient fields $K\subset K'$, respectively, such that $D'\cap K=D$. Let $p\in D[X]$ be monic and $d\in D$. Let $g\in D[X]$ be such that $g(X)$ is divisible in $D'[X]$ by $p(X)$ modulo $dD'[X]$. Then $g(X)$ is divisible in $D[X]$ by $p(X)$ modulo $dD[X]$.
\end{Lemma}

\vskip0.2cm

\begin{Rem} Equivalently, the lemma says that the ideal of $D[X]$ generated by $d$ and $p(X)$ is a contracted ideal with respect to  $D'[X]$, that is:
$$dD[X]+p(X)D[X]=(dD'[X]+p(X)D'[X])\cap D[X].$$
We also notice that if $d=0$ we have $p(X)D[X]=p(X)D'[X]\cap D[X]$ for every extension of domains $D\subset D'$. In this case we only use the  assumption that $p(X)$ is monic. 
\end{Rem}

\vskip0.3cm
\Proof As already observed in \cite{G}, the equality $D'\cap K=D$ is equivalent to the condition that every principal ideal of $D$ is contracted with respect to $D'$, that is $dD'\cap D=dD$, for every $d\in D$ (the proof is straightforward). From this fact it easily follows that $dD'[X]\cap D[X]=dD[X]$.\\
Let now $g\in D[X]$ be such that $g\in dD'[X]+p(X)D'[X]$. Then 
$$g(X)=Q(X)p(X)+dR(X)$$
for some $Q,R\in D'[X]$. By Euclidean division ($p(X)$ is monic) we can assume that $\deg(R)<\deg(p)$ (if not, just divide $R(X)$ by $p(X)$ in $D'[X]$). We divide now $g(X)$ by $p(X)$ in $D[X]$ and we get
$$g(X)=q(X)p(X)+r(X)$$
for some $q,r\in D[X]$, $\deg(r)<\deg(p)$. By the uniqueness of quotient and remainder in $D'[X]$ we have $Q(X)=q(X)$ and $dR(X)=r(X)$. This implies that $Q\in D[X]$. By the remark at the beginning of this proof, the second equality implies that $R\in D[X]$,  and this concludes the proof. $\Box$
\vskip0.5cm
\begin{Ex} In general the statement of Lemma \ref{gdpK} is false, if we drop the assumption $D'\cap K=D$. For example, let $D=\Z[\sqrt{5}]\subsetneq D'=O_K=\Z[(1+\sqrt{5})/2]\subset K=\Q(\sqrt{5})$, where $O_K$ is the ring of integers of the number field $K$ ($O_K$ is also equal to the integral closure of $D$ in $K$). Obviously, $D$ and $D'$ have the same quotient field $K$ and $D'\cap K\supsetneq D$. We have $2D\subsetneq 2O_K\cap D=2O_K\subset D$ (last equality holds because $2O_K$ is the conductor of the integral extension $D\subset O_K$). The element $\alpha=1+\sqrt{5}$ is in $2O_K$ but not in $2D$, so that the polynomial $g(X)=X+\alpha\in D[X]$ is divisible in $O_K[X]$ by $X$ modulo $2O_K[X]$ but it is not divisible in $D[X]$ by $X$ modulo $2D[X]$.
\end{Ex}

\begin{Rem} The Lemma \ref{gdpK} applies in the following settings: let $D$ be integrally closed, $E$ any  field extension of $K$ and $D_E$ the integral closure of $D$ in $E$. Then $D_E\cap K=D$, so the above statement holds with $D_E$ at the place of $D'$. Under these conditions, given a monic polynomial $p\in D[X]$ of degree $n$, Lemma \ref{crit0} and \ref{gdpK} together show that
$$\textnormal{Int}(M_n^{p}(D_E),M_n(D_E))\cap K[X]=\textnormal{Int}(M_n^{p}(D),M_n(D)).$$
In particular, if we extend the field $K$ to the splitting field $F$ of $p(X)$, then by Lemma \ref{gdpK} the same modular divisibility relation for a polynomial $g\in D[X]$ of condition 2) of Lemma \ref{crit0} holds. In this way, in the next section we will be able to apply the Newton expansion of a polynomial $g\in D[X]$ with respect to the set of roots of $p(X)$ in $D_F$ (the coefficients of this expansion are exactly the divided differences of $g(X)$ evaluated at the roots of $p(X)$). \end{Rem}
\vskip0.5cm

\section{Divided differences of integer-valued polynomials over matrices}\label{divdiff}

\vskip0.5cm
Let $g(X)$ be a polynomial of degree $n$ over a commutative ring $R$ and let $a_0,\ldots,a_n$ be arbitrary chosen elements of $R$ (possibly with repetitions). 

We have the following formula due to Newton:
\begin{align}\label{Newt}
g(X)=g(a_0)+\Phi^1(g)(a_0,a_1)(X-a_0)+\Phi^2(g)(a_0,a_1,a_2)(X-a_0)(X-a_1)+\ldots\nonumber\\
+\Phi^{n}(g)(a_0,\ldots,a_n)(X-a_0)\cdot\ldots\cdot(X-a_{n-1})
\end{align}
where, for each $k\in\N$, $\Phi^{k}(g)(X_0,\ldots,X_k)$ is the $k$-th divided difference of $g$ defined recursively as follows: 
\begin{align*}
\Phi^{0}(g)(X_0)\doteqdot& g(X)\\
\Phi^{1}(g)(X_0,X_1)\doteqdot& \frac{g(X_0)-g(X_1)}{X_0-X_1}\\
\ldots\\
\Phi^{k}(g)(X_0,\ldots,X_k)\doteqdot&\frac{\Phi^{k-1}(g)(X_0,\ldots,X_{k-1})-\Phi^{k-1}(g)(X_0,\ldots,X_{k-2},X_k)}{X_{k-1}-X_{k}}\\
\end{align*}
$\Phi^{k}(g)$ is a symmetric polynomial with coefficients in $R$ in $k+1$ variables. Obviously, a polynomial $g\in R[X]$ is divisible by $(X-a_0)\cdot\ldots\cdot(X-a_{n-1})$ if and only if the coefficients $\Phi^{k}(g)(a_0,\ldots,a_k)$ of the expansion in (\ref{Newt}) are equal to zero, for all $0\leq k<n$. We just use the fact that $\{1\}\cup\{\prod_{i=0,\ldots,k}(X-a_i)\}_{k=0,\ldots,n-1}$ are linearly independent over $D$. More in general, if $\{a_i\}_{i\in\N}$ is a sequence in $R$, then $\{1\}\cup\{\prod_{i=0,\ldots,k}(X-a_i)\}_{k\in\N}$ is a free basis of the $R$-module $R[X]$. 

The next lemma appears also in \cite[Proposition 11]{EvFaJoh}, but for the sake of the reader we report it here.

\begin{Lemma}[\cite{EvFaJoh}]\label{gdivpmodd}
Let $D$ be a domain. Let $g\in D[X]$, $d\in D\setminus\{0\}$ and $a_0,\ldots,a_n\in D$ (not necessarily distinct). Then $g(X)$ is divisible modulo $dD[X]$ by $p(X)=\prod_{i=0,\ldots,n}(X-a_i)$ if and only if for all $0\leq k<n$ we have $\Phi^k(g)(a_{0},\ldots,a_{k})\in dD$.
\end{Lemma}

\Proof Consider everything modulo $dD$ and apply Newton's formula. $\Box$
\vskip0.5cm
Let $D$ be an integrally closed domain with quotient field $K$. Remember that for any extension $E$ of $K$ we denote by $D_E$ the integral closure of $D$ in $E$. From now on we consider only algebraic extension of $K$, which are tacitly assumed to be contained in $\overline{K}$. We recall that a multi-set is a set with repetitions. For instance, the multi-set of roots of a polynomial $p\in K[X]$ is the multi-set formed by the roots of the polynomial
in $\overline{K}$ each root being repeated a number of times equals to its multiplicity in $p(X)$. We give now the following proposition:
\vskip0.6cm
\begin{Prop}\label{gdpdiv}
Let $D$ be an integrally closed domain with quotient field $K$. Let $g\in D[X]$, $d\in D\setminus\{0\}$ and $p\in D[X]$ monic of degree $n$. Let $\Omega_p=\{\alpha_1,\ldots,\alpha_n\}$ be the multi-set of roots of $p(X)$ in a splitting field $F$ over $K$. Then $f(X)=g(X)/d\in\textnormal{Int}(M_n^{p}(D),M_n(D))$ if and only if for all $0\leq k<n$, $\Phi^k(g)(\alpha_{i_0},\ldots,\alpha_{i_k})\in dD_F$, for all subsets $\{i_0,\ldots,i_k\}$ of $k+1$ elements of $\{1,\ldots,n\}$.
\end{Prop} 
\vskip0.2cm
\Proof Suppose $f(X)=g(X)/d$, for some $g\in D[X]$ and $d\in D$, $d\not=0$. By Lemma \ref{crit0}, $f(X)$ is in $\textnormal{Int}(M_n^{p}(D),M_n(D))$ if and only if $g(X)$ is divisible modulo $dD[X]$ by $p(X)$. Since $D$ is integrally closed, by Lemma \ref{gdpK} this is equivalent to $g(X)$ divisible modulo $dD_F[X]$ by $p(X)$. If we fix an order of the roots $\alpha_1,\ldots,\alpha_n$ of $p(X)$, by Newton's formula (\ref{Newt}) and Lemma \ref{gdivpmodd}, the previous condition holds if and only if $\Phi^k(g)(\alpha_{0},\ldots,\alpha_{k})\in dD_F$ for all $0\leq k<n$. Since this holds for every possible order of the $\alpha_i$'s we may choose (because the above modular divisibility condition does not depend on the order of the roots $\Omega_p$), we get the final statement. $\Box$


\vskip0.6cm
\begin{Rem}\label{Remgdpdiv} Notice that the roots of $p(X)$ are integral over $D$ (hence they are in $D_F$) and by definition $F=K(\alpha_1,\ldots,\alpha_n)$. In particular, for each of the relevant $k$, $\Phi^k(f)(\alpha_{i_0},\ldots,\alpha_{i_k})\in D_{F_k}$, where $F_k$ is the field generated by $\alpha_{i_0},\ldots,\alpha_{i_k}$ over $K$. We stress that we are not assuming $p(X)$ to be irreducible. Notice also that the above condition "$\Phi^k(g)(\alpha_{i_0},\ldots,\alpha_{i_k})\in dD_F$" is equivalent to "$\Phi^k(f)(\alpha_{i_0},\ldots,\alpha_{i_k})$ is integral over $D$" (since, obviously, $\Phi^k(f)(\alpha_{i_0},\ldots,\alpha_{i_k})$ belongs to $F$). If the conditions on the values of the divided differences $\Phi^k(g)$, for $0\leq k<n$, at the roots $\Omega_p$ of $p(X)$ hold for a fixed ordering, then they hold for every other order we may choose.
\end{Rem}
\vskip0.3cm
\begin{Ex}\label{Esn=2} Let us see with an example what happens in the case $n=2$. Let $p\in D[X]$ be an irreducible monic polynomial of degree $2$. Let $\alpha$ be a root of $p(X)$ in $\overline{K}$. Over the quadratic extension $K(\alpha)$ of $K$ we have  $p(X)=(X-\alpha)(X-\alpha')$, where $\alpha'$ is the conjugate root of $\alpha$. Then, given a polynomial $f(X)$ in $K[X]$, $f(X)$ is integer-valued over the set of matrices in $M_2(D)$ with characteristic polynomial equal to  $p(X)$ if and only if $f(\alpha)$ and $\Phi^1(f)(\alpha,\alpha')$ are integral over $D$ (by above, these two conditions implies that also $f(\alpha')$ is integral over $D$). If the extension $K(\alpha)\supset K$ is separable, then $\Phi^1(f)(\alpha,\alpha')\in D$, because $\Phi^1(f)(X,Y)$ is symmetric, so in particular $\Phi^1(f)(\alpha,\alpha')$ is invariant under the action of the Galois group of $K(\alpha)$ over $K$. If we consider a split polynomial $p(X)$ in $D[X]$, say $p(X)=(X-a)(X-b)$, for some $a,b\in D$, then $f\in\textnormal{Int}(M_2^{p}(D),M_2(D))$ if and only if $f(a)$ and $\Phi^1(f)(a,b)$ are in $D$.
\end{Ex}
\vskip0.6cm

Finally, by Proposition \ref{gdpdiv} and the representation of $\textnormal{Int}(M_n(D))$ in (\ref{IntMnIntMnpnozerJAc}) and (\ref{IntMnIntMnpzerJAc}), we give this characterization of the polynomials of $\textnormal{Int}(M_n(D))$ in terms of their divided differences.
\vskip0.5cm

\begin{Theorem}\label{MTINTMNGLOBALE}
Let $D$ be an integrally closed domain. Let $f\in K[X]$ and $n$ a positive integer. Then $f\in \textnormal{Int}(M_n(D))$ if and only if for every $0\leq k<n$, for every monic polynomial $p\in D[X]$ of degree $n$ and for every multi-set of $k+1$ elements $\{\alpha_1,\ldots,\alpha_{k+1}\}$ of the roots of $p(X)$ in $\overline K$, $\Phi^k(f)(\alpha_1,\ldots,\alpha_{k+1})$ is integral over $D$. 

\noindent If the intersection of the maximal ideals of finite index is $(0)$ then it is sufficient to consider polynomials $p\in D[X]$ which are irreducible.
\end{Theorem}
\Proof By representation (\ref{IntMnIntMnpnozerJAc}) and Proposition \ref{gdpdiv}, $f\in \textnormal{Int}(M_n(D))$ if and only if for all $0\leq k<n$, $\Phi^k(f)(\alpha_{i_0},\ldots,\alpha_{i_k})$ is integral over $D$ for every multi-set of $k+1$ elements of the roots $\Omega_p$ of any $p\in \mathcal P_n$. Notice that, by the arguments given in the proof of Proposition \ref{gdpdiv} and the Remark \ref{Remgdpdiv}, if for every $p\in \mathcal P_n$ we fix an ordering of the roots $\Omega_p$ of $p(X)$, for $f(X)$ to belong to $\textnormal{Int}(M_n^{p}(D),M_n(D))$ it is sufficient that each $k$-th divided difference of $f(X)$ is integral on the first $k+1$ roots of $p(X)$. The last claim follows by (\ref{IntMnIntMnpzerJAc}). $\Box$

\vskip0.6cm

\begin{Rem} The roots of any $p\in\mathcal{P}_n$ have degree over $K$ bounded by $n$. Notice that the divided differences $\Phi^k(f)$, for $0\leq k<n$, of an integer-valued polynomial $f(X)$ over $M_n(D)$ are integral on multi-sets $\{\alpha_1,\ldots,\alpha_{k+1}\}$ of integral elements over $D$ which are roots of polynomials in $\mathcal{P}_n$, that is, the elements $\{\alpha_1,\ldots,\alpha_{k+1}\}$ cannot be chosen independently (except when their degree is one, that is, they lie in $D$; this happens if $p(X)$ is totally split over $D$). Moreover, if  the intersection of the maximal ideals of finite index is $(0)$, we may just consider polynomials in $\mathcal P_n^{{\rm irr}}$, so that the $k$-th divided differences of $f(X)$ are integral on every set of $k+1$ conjugate integral elements of degree $n$ over $D$ inside $\overline{K}$.
\end{Rem}
\vskip0.3cm
\begin{Ex} Here we continue the Example \ref{Esn=2}. A polynomial $f\in K[X]$ is integer-valued over $M_2(D)$ if and only if the following two conditions hold:
\begin{itemize}
 \item[i)] for every $\alpha\in \overline{K}$, integral over $D$ of degree $2$, $f(\alpha)$ and $\Phi^1(f)(\alpha,\alpha')$ are integral over $D$ (where as before $\alpha'$ is the conjugate of $\alpha$ over $K$).
 \item[ii)] for every $(a,b)\in D^2$, $f(a),\Phi^1(f)(a,b)$ are in $D$.
\end{itemize}
If  the intersection of the maximal ideals of finite index is $(0)$, then $f\in{\rm Int}(M_2(D))$ if and only if condition i) holds.
\end{Ex}

\vskip0.5cm
\section{Integer-valued polynomials over triangular matrices with prescribed characteristic polynomial}

\vskip0.3cm

We give here some characterizations of integer-valued polynomials over triangular matrices. Now $D$ is again just an integral domain (that is, not necessarily integrally closed). The following definition is similar to Definition \ref{Mnp}.
\begin{Def}
Let $\mathcal{P}_n^s$ be the set of monic polynomials in $D[X]$ of degree $n$ which are totally split in $D$. Given a subset $P\subset\mathcal{P}_n^s$ we denote by 
$$T_n^P(D)\doteqdot\{T\in T_n(D)\,|\,p_T(X)\in P\}$$
the set of triangular matrices whose characteristic polynomial is in $P$. For a single polynomial $p\in\mathcal{P}_n^s$, we set $T_n^p(D)=T_n^{\{p\}}(D)$.
\end{Def}

\noindent The next lemma is analogous to \cite[Proposition 15]{EvFaJoh}.

\vskip0.4cm
\begin{Lemma}\label{crit0T}
 Let $p\in D[X]$ be a totally split polynomial of degree $n$. Let $f(X)=g(X)/d$ be in $K[X]$. Then
 $$f\in\textnormal{Int}(T_n^p(D),M_n(D))\Leftrightarrow g(X)\textnormal{ is divisible by }p(X)\textnormal{ modulo }dD[X].$$
In particular,
$$\textnormal{Int}(T_n^p(D),M_n(D))=\textnormal{Int}(M_n^p(D),M_n(D)).$$
\end{Lemma}

\Proof One direction is easy (the proof is exactly like in Lemma \ref{crit0}): suppose $g(X)=h(X)p(X)+dr(X)$, for some $h,r\in D[X]$. Let $T$ be a triangular matrix with characteristic polynomial equal to $p(X)$. Then by the Cayley-Hamilton theorem we have $g(T)=dr(T)$, so that $f(T)=r(T)\in M_n(D)$.

For the converse, we use an inductive argument on the degree of $p(X)$ as in \cite[Proposition 14]{EvFaJoh}. Let $n=2$ and let $p(X)=(X-a_0)(X-a_1)$, for some $a_0,a_1\in D$. Given $f(X)=g(X)/d$ in $\textnormal{Int}(T_2^p(D),M_2(D))$, we divide $g(X)$ by $p(X)$ in $D[X]$: $g(X)=q(X)p(X)+r(X)$, with $q,r\in D[X]$, $\deg(r)\leq 1$, say $r(X)=aX+b$. Take 
$$T=\left(\begin{matrix}
a_0 & 0 \\ 
1 & a_1
\end{matrix}\right).$$
Then $g(T)=aT+b\in T_2(dD)$ (essentially by (\ref{fMgN})) is the triangular matrix 
$$\left(\begin{matrix}
aa_0+b & 0 \\ 
a & aa_1+b
\end{matrix}\right)$$
so that $a,b\in dD$ as wanted.

Suppose now the statement is true up to $n$ and suppose that $p(X)=\prod_{i=0,\ldots,n}(X-a_i)$. Let $p_1(X)=\prod_{i=1,\ldots,n}(X-a_i)$. Given $T\in T_{n}^{p_1}(D)$, we consider the following $(n+1)\times (n+1)$ matrix:
$$T_1=\left(\begin{array}{c|c}
\raisebox{-15pt}{{\Large\mbox{{T}}}}& 0 \\ [-4ex] 
 & \vdots\\
 & 0 \\ \hline
0 \cdots 0 & a_0 
\end{array}\right).$$
Notice that $T_1$ is in $T_{n+1}^p(D)$ and it is a block matrix. Let $f\in\textnormal{Int}(T_{n+1}^p(D),M_{n+1}(D))$. Then $f(X)$ is integer-valued on $T$. In fact,
$$f(T_1)=\left(\begin{array}{c|c}
\raisebox{-10pt}{{\Large\mbox{{f(T)}}}}& 0 \\ [-4ex] 
 & \vdots\\
 & 0 \\ \hline
0 \cdots 0 & f(a_0) 
\end{array}\right)$$
$f(T_1)$ has entries in $D$, so that $f(T)$ has entries in $D$, too. Since this holds for every $T$ in $T_{n}^{p_1}(D)$ and for every $f\in\textnormal{Int}(T_{n+1}^p(D),M_{n+1}(D))$, we have $\textnormal{Int}(T_{n+1}^p(D),M_{n+1}(D))\subset \textnormal{Int}(T_{n}^{p_1}(D),M_n(D))$.

By induction hypothesis, given $f(X)=g(X)/d\in \textnormal{Int}(T_{n+1}^p(D),M_{n+1}(D))$, we have that  $g(X)$ is divisible by $p_1(X)$ modulo $dD[X]$. Then we continue as in the proof of Proposition 14 of \cite{EvFaJoh}, essentially going through the characterization of Lemma \ref{gdivpmodd} together with Corollary 13 of \cite{EvFaJoh}. 

The last claim follows from Lemma \ref{crit0}. $\Box$ 

\vskip0.5cm

\noindent The following proposition follows from the obvious fact $T_n^S(D)=\bigcup_{p\in S}T_n^p(D)$. This is analogous to Proposition \ref{criterio}.
\vskip0.5cm
\begin{Prop}\label{criteriotriang}
Let $P$ be a set of monic totally split polynomials in $D[X]$ of degree $n$.  Let $f(X)=g(X)/d\in K[X]$. Then
$$f\in\textnormal{Int}(T_n^P(D),M_n(D))\Leftrightarrow g(X) \textnormal{ divisible modulo }dD[X] \textnormal{ by all }p\in P.$$
In particular, 
 $$\textnormal{Int}(T_n^P(D),M_n(D))=\textnormal{Int}(M_n^P(D),M_n(D)).$$
\end{Prop}
\noindent The last statement follows from Proposition \ref{criterio}.
\vskip0.4cm

\section{Polynomial closure of set of matrices}\label{polclomatrices}
\vskip0.4cm
We recall from \cite{EvFaJoh} the following definition. Given a subset $S$ in $M_n(D)$, the polynomial closure of $S$ in $M_n(D)$ is the set
$$\overline{S}\doteqdot\{M\in M_n(D)\;|\;f(M)\in M_n(D),\;\forall f\in\textnormal{Int}(S,M_n(D))\}.$$
$\overline{S}$ is the largest set in $M_n(D)$ containing $S$ such that $\textnormal{Int}(S,M_n(D))=\textnormal{Int}(\overline{S},M_n(D))$. A set of matrices $S$ is said to be polynomially closed if $S=\overline{S}$, polynomially dense if $\overline{S}=M_n(D)$ (for example, $M_n^{\textnormal{irr}}(D)$ is polynomially dense in $M_n(D)$, if  the intersection of the maximal ideals of finite index is $(0)$, see (\ref{Mnirr})).

For some domains $D$, the authors of \cite{EvFaJoh} prove that the polynomial closure of $T_n(D)$ is equal to the set $S_n(D)$ of all matrices in $M_n(D)$ with totally split characteristic polynomial (see \cite[Theorem 27]{EvFaJoh}). On the other hand, Proposition \ref{criteriotriang} shows that, for any given subset $P\subseteq\mathcal{P}_n^s$, $T_n^P(D)$ and $M_n^P(D)$ have the same polynomial closure in $M_n(D)$, where $D$ is any integral domain.

Under a rather general assumption, we characterize the polynomial closure $\overline{M_n^p(D)}$ of $M_n^p(D)$ in $M_n(D)$, for any given monic polynomial $p\in D[X]$ of degree $n$. 
 
\begin{Prop}\label{polcloMnp}
 Suppose $D$ is an integral domain with zero Jacobson radical. Let $p\in \mathcal{P}_n$. Then
 $$\overline{M_n^p(D)}=\{M\in M_n(D)\,|\,p(M)=0\}.$$
Moreover, the eigenvalues (in $\overline{K}$) of every $M\in \overline{M_n^p(D)}$ are roots of $p(X)$.
 \end{Prop}
\Proof By Lemma \ref{crit0}, a polynomial $f(X)$ in $\textnormal{Int}(M_n^p(D),M_n(D))$ has the form $f(X)=r(X)+p(X)q(X)$, for some $q\in K[X]$, $r\in D[X]$, $\deg(r)<\deg(p)$. Clearly, if $p(M)=0$ then for every polynomial $f(X)$ as above, we have $f(M)=r(M)\in M_n(D)$. Conversely, a matrix $M$ is in the polynomial closure of $M_n^p(D)$ if and only if, for every $d\in D\setminus\{0\}$ and for every $g\in D[X]$, we have $p(M)g(M)\in M_n(dD)$. In particular (choose $g(X)=1$), this implies that $p(M)$ has entries which lie in $\cap_{d\in D\setminus\{0\}}dD$. Since $D$ has zero Jacobson radical, we have $p(M)=0$, as wanted. 

For the last claim, let $M$ be a matrix such tht $p(M)=0$. By a result of McCoy (see Corollary p. 283 of \cite{McC}) the ideals $N_{D[X]}(M)$ and $(p_M(X))$ of $D[X]$ have the same radical, so that $p(X)^m\in(p_M(X))$ for some $m\in\N$.  In particular, for all the roots $\alpha$ of $p_M(X)$ we have $p(\alpha)=0$. $\Box$ 
\vskip0.4cm
\noindent Without the above assumption on $D$, the polynomial closure of $M_n^p(D)$ is equal to the set $\{M\in M_n(D)\,|\,p(M)\in M_n(\mathcal{P})\}$, where $\mathcal{P}=\cap_{d\in D\setminus\{0\}}dD$.
\vskip0.5cm
\begin{Cor}
Let $D$ be an integral domain with zero Jacobson radical. Let $p\in \mathcal{P}_n$ be such that $p(X)$ is irreducible in $K[X]$. Then $M_n^p(D)$ is polynomially closed.
\end{Cor}
If Gauss's lemma applies (or just if $D$ is integrally closed, since $p(X)$ is monic), we may just suppose that $p(X)$ is irreducible in $D[X]$.
\vskip0.2cm
\Proof Let $M\in M_n(D)$ be a matrix in the polynomial closure of $M_n^p(D)$. By Lemma \ref{polcloMnp}, $p(M)=0$. Since $p(X)$ is irreducible it is equal to the minimal polynomial of $M$ (when consired in $M_n(K)$; a matrix over a commutative ring $R$ may not have a minimal polynomial). Then the characteristic polynomial of $M$ is equal to the minimal polynomial of $M$, because they always have the same irreducible factors. This means that $M$ is in $M_n^p(D)$. $\Box$
\vskip0.4cm
We notice that the above argument shows that $M_n^p(D)=\{M\in M_n(D)\,|\,p(M)=0\}$, under only the assumption that $p(X)$ is irreducible in $K[X]$ (thus, we don't need for this equality the fact that $D$ has zero Jacobson radical).

We leave the following open  problem.

\textbf{Problem}: given any subset $S$ in $\mathcal{P}_n$, compute $\overline{M_n^S(D)}$.
\vskip0.4cm
Since $M_n^S(D)=\bigcup_{p\in S}M_n^p(D)$, by Lemma \ref{polcloMnp} and known properties of polynomial closure (see \cite{CaCh}) we have
$$\bigcup_{p\in S}\overline{M_n^p(D)}=\{M\in M_n(D)\,|\,p(M)=0,\textnormal{ for some }p\in S\}\subset  \overline{M_n^S(D)}.$$

\begin{Ex} It is not true in general that the above containment is an equality. Let $D=\Z$ and let $K$ be a quadratic number field with ring of integers $O_K$. Let 
$M_2^K(\Z)$ be the set of matrices in $M_2(\Z)$ whose characteristic polynomial is the minimal polynomial of some $\alpha\in O_K\setminus\Z$. Then by Theorem 4.1 of \cite{Per} we have  $\textnormal{Int}(M_2^K(\Z),M_2(\Z))\subset\textnormal{Int}_{\Q}(O_K)\subset \IZ$, where $\textnormal{Int}_{\Q}(O_K)\doteqdot{\rm Int}(O_K)\cap \Q[X]$. In particular this implies that any $2\times 2$ diagonal  matrix $D$ with rational integers in the diagonal is contained in the polynomial closure of $M_2^K(\Z)$. Obviously, the null ideal of such a matrix $D$ over $\Z$ is generated by a totally split polynomial in $\Z$, which is not the minimal polynomial of some quadratic algebraic integer of $K$.
\end{Ex}
\vskip0.5cm

\section{Final Remarks}
\vskip0.5cm

Let $f\in{\rm Int}(M_n(D))$. As we anticipated in the introduction, for every $0\leq k<n$, $\Phi^k(f)$ is integer-valued on every multi-set $\{\alpha_0,\ldots,\alpha_k\}$ of elements of $D$: it is sufficient to consider polynomials which split into linear factors in $D[X]$ in Theorem \ref{MTINTMNGLOBALE}. This is exactly the reason why our result generalizes \cite[Theorem 16]{EvFaJoh}. This implies that $\textnormal{Int}(M_n(D))$ is a subring of the ring of integer-valued polynomials whose $k$-th divided differences are also integer-valued for all $k<n$. This ring has been introduced by Bharghava in \cite{BG} and it is denoted by:
$$\textnormal{Int}^{\{n-1\}}(D)\doteqdot\{f\in K[X]\,|\,\Phi^k(f)(D^{k+1})\subset D\,,\forall 0\leq k<n\}.$$
By \cite[Theorem 16]{EvFaJoh} the latter ring is equal to the ring of polynomials which are integer-valued over the set $T_n(D)$ of triangular matrices over $D$. 

By the same Theorem \ref{MTINTMNGLOBALE}, for $k=0$, we see that $\Phi^0(f)(X)=f(X)$ takes integral values on every element $\alpha$ which is integral over $D$ of degree $\leq n$. In particular, $f\in\textnormal{Int}_{K}(D_E)\doteqdot K[X]\cap\textnormal{Int}(D_E)$, for every finite extension $E$ of $K$ of degree $\leq n$ (as usual, $\textnormal{Int}(D_E)$ is the subring of $E[X]$ of those polynomials $f(X)$ such that $f(D_E)\subset D_E$). In this way we have the following containment: 
$$\textnormal{Int}(M_n(D))\subset\bigcap_{[E:K]\leq n}\textnormal{Int}_{K}(D_E).$$
In the case $D=\Z$ this inclusion was showed in \cite{LopWer}. It can be easily proved that the latter intersection is equal to the ring $\textnormal{Int}_{K}(\mathcal{A}_n)\doteqdot\{f\in K[X]\,|\,f(\mathcal{A}_n)\subset\mathcal{A}_n \}$, where $\mathcal{A}_n$ is the set of all the integral elements over $D$ in $\overline{K}$ of degree bounded by $n$ (see also \cite{Per}, in the case $D=\Z$).
\vskip0.7cm
\subsection*{\textbf{Acknowledgments}}
I wish to thank Keith Johnson for the useful suggestions. I also thank the referee for the several suggestions he/she gave which improved the overall quality of the paper. 
The author was supported by the Austrian Science Foundation (FWF), Project Number P23245-N18.

\end{document}